\documentclass[BCOR8mm]{scrartcl}
\usepackage{epsfig} 
\usepackage{amsmath} 
\usepackage{amsthm}
\usepackage{fancyhdr}
\usepackage{amssymb}  
\usepackage[mathscr]{euscript}
\usepackage{nicefrac}
\newcommand{\field}[1]{\mathbb{#1}}
\newcommand{\R}{\field{R}}
\newtheorem{theorem}{Theorem}
\newtheorem{corollary}[theorem]{Corollary}
\newtheorem{lemma}[theorem]{Lemma}
\newtheorem{proposition}[theorem]{Proposition}
\newtheorem{remark}[theorem]{Remark}
\newtheorem{definition}[theorem]{Definition}
\newtheorem{example}[theorem]{Example}
\newcommand{\einsnorm}[2]{\ensuremath{
    \!\!\;\!\!\!\;
    \left\bracevert\!\!\!\!\!\left\bracevert
    \!
        #1(#2)
    \!
      \right\bracevert\!\!\!\!\!\right\bracevert
    \!\!\;\!\!\!\;
  }}
\newcommand{\supnorm}[2]{\ensuremath{
    \!\!\;\!\!\!\;
    \left\bracevert\!\!\!\!\left\bracevert
        {#1}_{#2}
      \right\bracevert\!\!\!\!\right\bracevert
    \!\!\;\!\!\!\;
  }}

\newcommand{\einsbisn}{\{1,\ldots,n\}}
\let\epsilon=\varepsilon

\begin{document}
\title{An ISS Small-Gain Theorem for General Networks}
\author{ Sergey Dashkovskiy\thanks{
    $,{}^\dag$
    \-Universit\"at \-Bremen, \-Zentrum \-f\"ur \-Technomathematik, 
    Postfach 330440,
    28334 Bremen, \-Germany,
    Phone +49-421-218-9407, Fax +49-421-218-4235,
    \texttt{\{dsn,rueffer\}@math.uni-bremen.de}}
  \and Bj\"orn S. R\"uffer\footnotemark{} 
  \and 
  Fabian R.  Wirth\thanks{
    The Hamilton Institute, NUI Maynooth, Maynooth, Co.~Kildare\-,
    Ireland, \texttt{fabian.wirth@nuim.ie}
  }
}

\date{June 7, 2005}

\maketitle

\begin{abstract}
  We provide a generalized version of the \emph{nonlinear small-gain
    theorem} for the case of more than two coupled input-to-state
  stable (ISS) systems. For this result the interconnection gains are
  described in a nonlinear gain matrix and the small-gain condition
  requires bounds on the image of this gain matrix. The condition may
  be interpreted as a nonlinear generalization of the requirement that
  the spectral radius of the gain matrix is less than one.  We give
  some interpretations of the condition in special cases covering two
  subsystems, linear gains, linear systems and an associated
  artificial dynamical system.
\end{abstract}

\paragraph{Keywords}
{\newcommand{\and}{~--~}  
 Interconnected systems \and input-to-state stability \and
  small-gain theorem \and large-scale systems \and monotone maps }

\paragraph{MSC-classification:} 93C10 (Primary) 34D05, 90B10, 93D09, 93D30 (Secondary)

\section{Introduction}
%auto-ignore
Stability is one of the fundamental concepts in the analysis and design of
nonlinear dynamical systems. The notions of input-to-state stability (ISS) and
nonlinear gains have proved to be an efficient tool for the qualitative
description of stability of nonlinear input systems. 
There are different equivalent formulations of ISS: In terms
of  $\mathscr{KL}$ and $\mathscr K_\infty$ functions (see below),
via Lyapunov functions, as an asymptotic stability property combined with
asymptotic gains, and others, see \cite{SoW96}. A more
quantitative but equivalent formulation, which captures the long
term dynamic behavior of the system, is the notion of
input-to-state dynamical stability (ISDS), see \cite{Gru02}. 

One of the interesting properties in the study of ISS systems is that under
certain conditions input-to-state stability is preserved if ISS systems are
connected in cascades or feedback loops.  In this paper we generalize the
existing results in this area. In particular, we obtain a general condition
that guarantees input-to-state stability of a general system described as an
interconnection of several ISS subsystems.

The earliest interconnection result on ISS systems states that
cascades of ISS systems are again ISS, see e.g.,
\cite{SoT95,Son89,Son01}.  Furthermore, small-gain theorems for the
case of two ISS systems in a feedback interconnection have been
obtained in
\cite{Gru02,JTP94,JMW96}. These results state in one way or another that if
the composition of the gain functions of ISS subsystems is smaller than the
identity, then the whole system is ISS. 

The papers \cite{Gru02,JTP94,JMW96} use different approaches to the
formulation of small-gain conditions that yield sufficient stability criteria:
In \cite{JTP94} the proof is based on the properties of $\mathscr{KL}$ and
$\mathscr K_\infty$ functions. This approach requires that the composition of
the gains is smaller than the identity in a robust sense, see below for the
precise statement. We show in Example~\ref{JTP-example} that within the
context of this approach the robustness condition cannot be
weakened.  The result in that paper also covers practical ISS results, which
we do not treat here. An ISS-Lyapunov function for the feedback system is
constructed in \cite{JMW96} as some combination of the corresponding
ISS-Lyapunov functions of both subsystems.  The key assumption of the proof in
that paper is that the gains are already provided in terms of the Lyapunov
functions, by which the authors need not resort to a robust version of the
small-gain condition.  The proof of the small-gain theorem in \cite{Gru02} is
based on the ISDS property and conditions for asymptotic stability of the
feedback loop without inputs are derived.
These results will turn out to be special cases of
our main result. 

General stability conditions for large scale interconnected systems have been
obtained by various authors in other contexts.  In \cite{MiM77} sufficient
conditions for the asymptotic stability of a composite system are stated in
terms of the negative definiteness of some test matrix. This matrix is defined
through the given Lyapunov functions of the interconnected subsystems.
Similarly, in \cite{RoucHabe77} conditions for the stability of interconnected
systems in terms of Lyapunov functions of the individual systems are obtained.

In \cite{Sil79} {\v{S}}iljak considers structural perturbations and
their effects on the stability of composite systems using Lyapunov
theory.  The method is to reduce each subsystem to a one-dimensional
one, such that the stability properties of the reduced aggregate
representation imply the same stability properties of the original
aggregate system. In some cases the aggregate representation gives
rise to an interconnection matrix $\bar W$, such that quasi dominance
or negative definiteness of $\bar W$ yield asymptotic stability of the
composite system.

In \cite{Vid81} small-gain type theorems for general interconnected systems
with linear gains can be found.  These results are of the form that the
spectral radius of a gain matrix should be less than one to conclude
stability.  The result obtained here may be regarded as a nonlinear
generalization in the same spirit.

In this paper we consider a system which consists of two or more ISS
subsystems. We provide conditions by which the stability question of the
overall system can be reduced to consideration of stability of the subsystems.
We choose an approach using estimates involving $\mathscr{KL}$ and $\mathscr
K_\infty$ functions to prove the ISS stability result for general interconnected
systems. The generalized small-gain condition we obtain is, that for some
monotone operator $\tilde\Gamma$ related to the gains of the individual
systems the condition
\begin{equation}
  \label{eq:intro2}
  \tilde\Gamma (s) \not\geq s
\end{equation}
holds for all $s\geq0, s\ne 0$ (in the sense of the component-wise ordering of
the
positive orthant). 
We discuss interpretations of this condition in
Section~\ref{sec:interpretations}.

Although we believe our approach to be amenable to the
explicit construction of a Lyapunov function given the ISS-Lyapunov
functions for the subsystems, so far we have been able to prove this
only for linear gains.

While the general problem can be approached by repeated application of
the cascade property and the known small-gain theorem, in general this
can be cumbersome and it is by no means obvious in which order
subsystems have to be chosen to proceed in such an iterative manner.
Hence an extension of the known small-gain theorem
to larger interconnections is needed.

In this paper we obtain this extension for the general case.
Further, we show how to calculate the
gain matrix for linear systems and give some interpretation of
our result.

The paper is organized as follows. In
Section~\ref{sec:motiv-probl-descr} we introduce notation and
necessary concepts and state the problem. In particular, we will need
some basic properties of the positive orthant $\R_+^n$ interpreted as
a lattice.  In Section 3 we prove the main result, which generalizes
the known small-gain theorem, and consider the special case of linear
gains, for which we also construct an ISS-Lyapunov function. In
Section~\ref{sec:interpretations} the small-gain condition of the main
result is discussed and we show in which way it may be interpreted as
an extension of the linear condition that the spectral radius of the
gain matrix has to be less than one. There we also point out the
connection to some induced monotone dynamical system.  In
Section~\ref{sec:appl-line-syst} we show how the gain matrix can be
found for linear systems.  We conclude with
Section~\ref{sec:conclusions}.

\section{Problem description}
\label{sec:motiv-probl-descr}
%auto-ignore
\renewcommand{\labelenumi}{(\roman{enumi})}

\paragraph{Notation}
By $x^T$ we denote the transpose of a vector $x\in\R^n $.
For $x,y\in\R^n$,  we use the following notation
\begin{gather}
  x\ge y\;\Leftrightarrow \; x_i\ge y_i, \; i=1,\dots,n\,,\text{ and }
  x> y\;\Leftrightarrow \; x_i> y_i, \; i=1,\dots,n.  
\end{gather}
In the following $\R_+:=[0,\infty)$ and by $\R^n_+$ we denote $\{x\in\R^n:
x\geq 0\}$.  For a function $v:\R_+\to\R^m$ we define its restriction to the
interval $[s_1,s_2]$ by
$$v_{[s_1,s_2]}(t):=
\begin{cases}
  v(t) & \mbox{~if~}t\in[s_1,s_2],\\
  0 & \mbox{~else.}
\end{cases} $$
\begin{definition}
  (i)
  A function $\gamma:\mathbb{R}_+\to\mathbb{R}_+$ is said to
    be of class $\mathscr{K}$ if it is continuous, increasing and
    $\gamma(0)=0$. It is of class $\mathscr{K}_\infty$ if, in
    addition, it is proper, i.e., unbounded. 
  
    (ii)
  A function
    $\beta:\mathbb{R}_+\times\mathbb{R}_+\to\mathbb{R}_+$ is said
    to be of class $\mathscr{KL}$ if, for each fixed $t$, the function
    $\beta(\cdot,t)$ is of class $\mathscr{K}$ and, for each fixed
    $s$, the function $\beta(s,\cdot)$ is non-increasing and tends to
    zero for $t\to \infty$.
\end{definition}
Let $|\cdot|$ denote some norm in $\R^n$, and let in particular $|x|_{\max}=
\max_i |x_i|$ be the maximum norm. The
essential supremum norm on essentially bounded functions defined on $\R_+$ is
denoted by $\Vert\cdot\Vert_\infty$.
\begin{definition}\label{def:iss}
Consider a system
$$
\dot x=f(x,u),\;x\in\R^n, u\in\R^m $$
such that for all initial values
$x_0$ and all essentially bounded inputs $u$ unique solutions exist for all
positive times. We denote these solutions by $\xi(t;x_0,u)$.  The system is
called input to state stable (ISS), if there exist functions $\beta$ of class
${\mathscr {KL}}$ and $\gamma$ of class ${\mathscr K}$, such that the
inequality 
$$  |\xi(t;x_0,u)|\leq \beta(|x_0|,t)+\gamma(||u||_\infty) $$
holds for all $t\geq0, x_0 \in \R^n, u: \R_+\to \R^m$ essentially bounded.
\end{definition}

%auto-ignore
\paragraph{Problem statement}
Consider $n$ interconnected control systems given by
\begin{equation}
  \label{eq:ps1}
  \begin{array}{c}\dot x_1=f_1(x_1,\ldots,x_n,u)\\
  \vdots\\
  \dot x_n=f_n(x_1,\ldots,x_n,u)
  \end{array}
\end{equation}
where $x_i\in\R^{N_i}, u\in\R^{L}$ and
$f_i:\R^{\sum_{j=1}^nN_j+L}\to\R^{N_i}$ is continuous and
Lipschitz in the first $n$ arguments uniformly with respect to $u$
for $i=1,\ldots,n$. Here $x_i$ is the state of the
$i^{\mbox{\small th}}$ subsystem, and $u$ is considered as an
external control variable. 

We may consider $u$ as partitioned $u=(u_1,\ldots,u_n)$, such that
each $u_i$ is the input for subsystem $i$ only. Then each $f_i$ is of
the form $f_i(\ldots,u)=\tilde f_i(\ldots,P_i(u))=\tilde
f_i(\ldots,u_i)$ with some projection $P_i$. So without loss of
generality we may assume to have the same input for all systems.

We call the $i^{\mbox{th}}$ subsystem of (\ref{eq:ps1}) ISS, if there exist
functions $\beta_i$ of class ${\mathscr {KL}}$ and $\gamma_{ij},\gamma$ of
class ${\mathscr K}$, such that the solution $x_i(t)$ starting at
$x_i(0)$ satisfies
\begin{equation}
  \label{eq:ps2}
  |x_i(t)|\leq \beta_i(|x_i(0)|,t)+\sum_{j=1}^n
   \gamma_{ij}(||{x_j}_{[0,t]}||_\infty)
   +\gamma(||u||_\infty)
\end{equation}
for all $t\geq0$.

For notational simplicity we allow the case $\gamma_{ij}\equiv 0$ and
require $\gamma_{ii}\equiv 0$ for all $i$.  The functions
$\gamma_{ij}$ and $\gamma$ are called (nonlinear) gains.  We define
$\Gamma:\R^n_+\to\R^n_+$ by
\begin{equation}
\Gamma:=(\gamma_{ij}), \qquad
\Gamma(s_1,\ldots,s_n)^T := \left(
  \sum_{j=1}^n \gamma_{1j}(s_j),\ldots,\sum_{j=1}^n \gamma_{nj}(s_j)
\right)^T\label{eq:ps_gamma}
\end{equation}
for $s=(s_1,\ldots,s_n)^T\in\R^n_+$. We refer to $\Gamma$ as the \emph{gain
  matrix}, noting that it does not represent a linear map.  Note that by the
properties of $\gamma_{ij}$ for $s_1,s_2\in\R_+^n$ we have the implication
\begin{equation}\label{eq:monot}
s_1\ge s_2\;\Rightarrow\;\Gamma(s_1)\ge\Gamma(s_2),
\end{equation}
so that $\Gamma$ defines a monotone map.\\
Assuming each of the subsystems of \eqref{eq:ps1} to be ISS, we are interested
in conditions guaranteeing that the whole system defined by
$x=(x^T_1,\ldots,x^T_n)^T, f=(f^T_1,\ldots,f^T_n)^T$ and
\begin{equation}
  \label{eq:ps3}
  \dot x=f(x,u)
\end{equation}
is ISS (from $u$ to $x$).

%auto-ignore
\paragraph{Additional Preliminaries}

We also need some notation from lattice theory, cf. \cite{Sza63} for
example.  Although $(\R^n_+,\sup,\inf)$ is a lattice, with $\inf$
denoting infimum and $\sup$ denoting supremum, it is not complete.
But still one can define the upper limit for bounded functions
$s:\R_+\to\R^n_+$ by $$
\limsup_{t\to\infty}s(t):= \inf_{t\geq0}
\sup_{\tau\geq t} s(\tau).$$

For vector functions
$x=(x_1^T,\ldots,x_n^T)^T:\R_+\to\R^{N_1+\ldots+N_n}$ such that\linebreak[3]
$x_i:\R_+\to\R^{N_i}, i=1,\ldots,n$ and times $0\leq t_1 \leq t_2$ we
define
\[
\supnorm{x}{[t_1,t_2]}
:=
  \begin{pmatrix}
    \|x_{1,[t_1,t_2]}\|_\infty\\
    \vdots\\
    \|x_{n,[t_1,t_2]}\|_\infty
  \end{pmatrix}\in\R^n_+
  .
\]
We will need the following property.

\begin{lemma}
  \label{limsuplemma} 
  Let $s:\R_+\to\R^n_+$ be continuous and bounded. Then
  (setting $N_i\equiv 1$)
  $$\limsup_{t\to\infty} s(t)=\limsup_{t\to\infty} \supnorm
  {s}{[\nicefrac t2,\infty)}.$$
\end{lemma}

  \begin{proof} Let $\limsup_{t\to\infty} s(t)=:a\in\R^n_+$
    and $\limsup_{t\to\infty} \supnorm {s}{[\nicefrac
      t2,\infty)}=:b\in\R^n_+$.  For every $\epsilon\in\R^n_+,$ \mbox{$\epsilon>0$}
    (component-wise!) there exist $t_a,t_b\geq0$ such that
    \begin{eqnarray}
      \label{eq:limsup a}
      \forall\, t\geq t_a: \sup_{t\geq t_a} s(t) \leq a+\epsilon
      \quad\mbox{and}\quad
      \forall\, t\geq t_b:  \sup_{t\geq t_b} \supnorm {s}{[\nicefrac
        t2,\infty)}\leq b+\epsilon.
    \end{eqnarray}
    Clearly we have $$s(t) \leq \supnorm {s}{[\nicefrac
      t2,\infty)}$$
    for all $t\geq 0$, i.e., $a\leq b$. On the other hand
    $s(\tau)\leq a+\epsilon$ for $\tau\geq t$ implies
    $\supnorm {s}{[\nicefrac \tau2,\infty)} \leq a+\epsilon$ for
    $\tau\geq2t$, i.e., $b\leq a$. This
    immediately gives $a=b$, and the claim is proved.
  \end{proof}

Before we introduce the ISS criterion for interconnected systems let us
briefly discuss an equivalent formulation of ISS. A system
\begin{equation}
  \label{eq:prelim-system}
  \dot x=f(x,u),
\end{equation}
with $f:\R^{N+L}\to\R^{N}$ continuous and Lipschitz in $x\in\R^N$,
uniformly with respect to $u\in\R^L$, is said to have the \emph{asymptotic
  gain property} (AG), if there exists a function $\gamma_{AG}\in\mathscr
K_\infty$ such that for all initial values
$x_0\in\R^{N}$ and all essentially bounded control functions $u(\cdot):\R_+\to\R^L$,
\begin{equation}
  \label{eq:AG}
  \limsup_{t\to 0}|x(t;x_0,u)|\leq\gamma_{AG}(\|u\|_\infty).
\end{equation}
The asymptotic gain property states, that every trajectory must ultimately
stay not far from zero, depending on the magnitude of $\|u\|_\infty$.

The system (\ref{eq:prelim-system}) is said to be \emph{globally asymptotically
  stable at zero} (0-GAS), if there exists a $\beta_{GAS}\in\mathscr {KL}$,
such that for all initial conditions $x_0\in\R^{N}$
\begin{equation}
  \label{eq:0-GAS}
  |x(t;x_0,0)|\leq\beta_{GAS}(|x_0|,t).
\end{equation}
Thus 0-GAS holds, if, when the input $u$ is set to zero,
the system (\ref{eq:prelim-system}) is globally asymptotically stable at
$x^*=0$.

By a result of Sontag and Wang \cite{SoW96} the asymptotic gain property and
global asymptotic stability at 0 together are equivalent to ISS.

\section{Main results}
%auto-ignore
\let\phi=\varphi

In the following subsection we present a nonlinear version of the
small-gain theorem for networks. In
Subsection~\ref{subsec:linear-gains-and-iss-lyapunov-function} we
restate this theorem for the case when the gains are linear functions.
Here we also provide a method on how to construct an ISS-Lyapunov
function for the whole network system from given ISS-Lyapunov
functions of the subsystems.

\subsection{Nonlinear gains}

We introduce the following notation. For $\alpha_i\in\mathscr
K_\infty, i=1,\ldots,n$ define $D:\R_+^n\to\R_+^n$ by
\begin{equation}
  \label{eq:definition von D}
  D (s_1,\ldots,s_n)^T :=
  \left(
    \begin{array}{c}
      (\mbox{Id}+\alpha_1)(s_1)\\
      \vdots\\
      (\mbox{Id}+\alpha_n)(s_n)
    \end{array}
  \right).
\end{equation}

\begin{theorem}[small-gain theorem for networks]
  \label{thm:small-gain_for_networks}
  Consider the system \eqref{eq:ps1} and suppose that each subsystem is ISS,
  i.e., condition (\ref{eq:ps2}) holds for all $i=1,\ldots,n$.  Let $\Gamma$
  be given by (\ref{eq:ps_gamma}).  If there exists a mapping $D$ as in
  (\ref{eq:definition von D}), such that
  \begin{equation}
    (\Gamma\circ D)(s)\not\geq s,\qquad\forall s\in\R^n_+\setminus\{0\}\,,
  \label{eq:Gamma D Id}
  \end{equation}
  then the system (\ref{eq:ps3}) is ISS from $u$ to $x$.
\end{theorem}

\begin{remark}
  Although looking very complicated to handle at first sight,
  condition (\ref{eq:Gamma D Id}) is a straightforward extension of the ISS
  small-gain theorem of \cite{JTP94}. It has many interesting
  interpretations, as we will discuss in Section
  \ref{sec:interpretations}.
\end{remark}

The following lemma provides an essential argument in the proof of
Theorem~\ref{thm:small-gain_for_networks}.
\begin{lemma}
  \label{lemma:boundedness_ID-Gamma}
  Let $D$ be as in (\ref{eq:definition von D}) and suppose (\ref{eq:Gamma D
    Id}) holds. Then there exists a $\phi\in \mathscr K_\infty$ such that for
  all $w,v \in \R^n_+$,
  \begin{equation}
    \label{eq:ineqbound}
    (\mbox{Id}-\Gamma)(w)\leq v
  \end{equation}
  implies $|w|\leq \phi(|v|)$.
\end{lemma}

\begin{proof}
  Fix $v \in \R^n_+$. We first show, that for those $w\in \R^n_+$ satisfying
  \eqref{eq:ineqbound} at least some components have to be bounded. To this
  end let
  \begin{equation}
    \label{eq:newr}
    \begin{aligned}
      r^* := (D - \mbox{Id})^{-1}(v) = \left(
      \begin{array}{c}
        \alpha_1^{-1}(v_1)\\
        \vdots\\
        \alpha_n^{-1}(v_n)
      \end{array}
    \right)  \\
    \text{ and }
    s^* := D(r^*)= \left(
      \begin{array}{c}
        v_1+\alpha_1^{-1}(v_1)\\
        \vdots\\
        v_n+\alpha_n^{-1}(v_n)
      \end{array}
    \right)\,.
  \end{aligned}
\end{equation}
We claim that $s\geq s^*$ implies that $w=s$ does not satisfy
\eqref{eq:ineqbound}.  So let $s\geq s^*$ be arbitrary and $r=D^{-1}(s)\geq
r^*$ (as $D^{-1} \in {\cal K}_\infty^n$). For such $s$ we have
  \[ s - D^{-1}(s) = D(r) - r \geq D(r^*) - r^* = v \,,\]
  where we have used that $(D-\mbox{Id}) \in {\cal K}_\infty^n$.
  The assumption that $w=s$ satisfies \eqref{eq:ineqbound} leads to
  \[ s \leq v + \Gamma(s) \leq s - D^{-1}(s) + \Gamma(s)\,,\]
  or equivalently, $0 \leq \Gamma(s) - D^{-1}(s)$. This implies
  for $r= D^{-1}(s)$ that
  \[ r \leq \Gamma \circ D(r) \,,\]
  in contradiction to (\ref{eq:Gamma D Id}). This shows that the set of $w\in
  \R^n_+$ satisfying \eqref{eq:ineqbound} does not intersect the set
  \[ Z_1 := \{ w \in \R^n_+ \;|\; w \geq s^* \} \,.\]

  Assume now that $w\in \R^n_+$ satisfies \eqref{eq:ineqbound}. Let
  $s^1:=s^*$.
  If $s^1\not\geq w$, then there exists an index set $I_1\subset\einsbisn$,
  such that
  \begin{eqnarray*}
    w_i>s_i^1,
    ~
    \mbox{for}
    ~
    i\in I_1\quad\mbox{and}
    \quad
    w_i\leq s_i^1, &&
    ~\mbox{for}
    ~
    i\in I_1^c:=\einsbisn\setminus I_1\,.
  \end{eqnarray*}

  For index sets $I$ and $J$ denote by $y_I$ the restriction
  $$y_I:=(y_i)_{i\in I}$$
  for vectors $y\in\R^n_+$ and by
  $A_{IJ}:\R^{\#I}_+\to\R^{\#J}_+$ the restriction
  $$A_{IJ}:=(a_{ij})_{i\in I,j\in J}$$
  for mappings
  $A=(a_{ij})_{i,j\in\einsbisn}:\R^n_+\to\R^n_+$.

  So from (\ref{eq:ineqbound}) we obtain
  \[
  \begin{bmatrix}
    w_{I_1} \\ w_{I_1^c}
  \end{bmatrix} -
  \begin{bmatrix}
    \Gamma_{{I_1}{I_1}} & \Gamma_{{I_1}{I_1^c}} \\
    \Gamma_{{I_1^c}{I_1}} & \Gamma_{{I_1^c}{I_1^c}}
  \end{bmatrix}
  \left(
  \begin{bmatrix}
    w_{I_1} \\ w_{I_1^c}
  \end{bmatrix}
  \right)
  \leq
\begin{bmatrix}
    v_{I_1} \\ v_{I_1^c}
  \end{bmatrix}
  \,.\]
  Hence we have in particular
  \begin{equation}
    \begin{split}
      w_{I_1}-&\Gamma_{I_1I_1}(w_{I_1}) \leq
      v_{I_1}+\Gamma_{I_1I_1^c}(s^1_{I_1^c})
      \\
      &\leq \underbrace{D_{I_1}\circ(D_{I_1}-\mbox{Id}_{I_1})^{-1}}_{>\mbox{Id}} \circ
      (v_{I_1}+\Gamma_{I_1I_1^c}(s^1_{I_1^c}))
      =: s^2_{I_1}.
    \end{split}
    \label{eq:s^2_{I_1}}
  \end{equation}
  Note that $\Gamma_{I_1I_1}$ satisfies \eqref{eq:Gamma D Id} with $D$
  replaced by $D_{I_1}$. Thus,
  arguing just as before, we obtain, that $w_{I_1} \geq s^2_{I_1}$ is not
  possible. Hence some more components of $w$ must be bounded.

  We proceed inductively, defining $$
  I_{j+1}\subsetneqq I_{j},\quad
  I_{j+1}:=\{i\in I_j:w_i>s_i^{j+1}\},$$
  with $I_{j+1}^c:=\einsbisn\setminus
  I_{j+1}$ and $$
  s^{j+1}_{I_j}:=
  D_{I_j}\circ(D_{I_j}-\mbox{Id}_{I_j})^{-1} \circ
  (v_{I_j}+\Gamma_{I_jI_j^c}(s^j_{I_j^c})).
  $$

  Obviously this nesting will end after at most $n-1$ steps: There
  exists a maximal $k\leq n$, such that $$\einsbisn \supsetneqq I_1
  \supsetneqq \ldots \supsetneqq I_{k} \ne \emptyset$$
  and all
  components of $w_{I_{k}}$ are bounded by the corresponding components
  of~$s^{k+1}_{I_{k}}$.
  For $i=1,\ldots,n$ define $$
  \zeta_i := \max \{ j\in\einsbisn: i\in I_j \}
  $$ and
  $$
  s_\zeta := (s_1^{\zeta_1},\ldots,s_n^{\zeta_n}).
  $$
  Clearly we have $$
  w\leq s_\zeta \leq [ D\circ(D-\mbox{Id})^{-1} \circ
  (\mbox{Id}+\Gamma) ] ^n (v)
  $$
  and the term on the very right hand side does not depend on any
  particular choice of nesting of the index sets. Hence every $w$ satisfying 
  \eqref{eq:ineqbound} also satisfies
  $$
  w \leq [ D\circ(D-\mbox{Id})^{-1} \circ
  (\mbox{Id}+\Gamma) ] ^n \circ
  \begin{pmatrix}
    |v|_{\max},&
    \ldots,&
    |v|_{\max}
  \end{pmatrix}^T
  $$
  and taking the $\max$-norm on both sides yields $$
  |w|_{\max}
  \leq \phi( |v|_{\max})$$
  for some function $\phi$ of class
  $\mathscr{K}_\infty$. This completes the proof of the lemma.
\end{proof}

We proceed with the proof of Theorem
\ref{thm:small-gain_for_networks}, which is divided into two main
steps. First we establish the existence of a solution of the system
(\ref{eq:ps3}) for all times $t\geq 0$. In the second step we
establish the ISS property for this system.

\begin{proof}(of Theorem~\ref{thm:small-gain_for_networks})
  \emph{Existence of a solution for (\ref{eq:ps3}) for all times:} For
  finite times 
  $t\geq 0$    
  and for
  $s\in\R^n_+$ we introduce the abbreviating notation
  \begin{gather}
    \einsnorm{x}{t}
    :=
    \left(
      \begin{array}{c}
        |x_1(t)|\\
        \vdots\\
        |x_n(t)|
      \end{array}
    \right)\in\R^n_+\,,
    \quad
   \gamma^n(\|u\|_\infty)
    :=
    \left(
      \begin{array}{c}
        \gamma(\|u\|_\infty) \\
        \vdots\\
        \gamma(\|u\|_\infty)
      \end{array}
    \right)\in\R^n_+\\
    \label{eq:definition beta}
    \mbox{and}\quad
    \beta(s,t)
    :=
    \left(
      \begin{array}{c}
        \beta_1(s_1,t)\\
        \vdots\\
        \beta_n(s_n,t)
      \end{array}
    \right):\R^n_+\times\R_+\to\R^n_+\,.
  \end{gather}
  Now we can rewrite the ISS conditions (\ref{eq:ps2}) of the
  subsystems in a vectorized form for $\tau\geq 0$ as
  \begin{equation}
    \label{eq:ISS subsysteme}
    \einsnorm{x}{\tau}\leq \beta(\einsnorm{x}{0},\tau) + \Gamma(
    \supnorm{x}{[0,\tau]})
    + \gamma^n(\|u\|_\infty)
  \end{equation}
  and taking the supremum on both sides over $\tau\in[0,t]$ we obtain
  \begin{equation}
    \label{eq:Subsysteme beschränkt}
    \begin{aligned}
      (\mbox{Id}-\Gamma)\circ\supnorm{x}{[0,t]}
      &=&
      \supnorm{x}{[0,t]}-\Gamma(\supnorm{x}{[0,t]})\\
      &\leq&
      \beta(\einsnorm{x}{0},0) + \gamma^n(\|u\|_\infty)
    \end{aligned}
  \end{equation}
  where we used (\ref{eq:monot}).
  Now by Lemma \ref{lemma:boundedness_ID-Gamma} we find
  \begin{equation}
    \|x_{[0,t]}\|_\infty\leq \phi\left(\left|\beta(\einsnorm{x}{0},0) +
    \gamma^n(\|u\|_\infty)\right|\right) =: s_\infty \label{eq:s_infty}
  \end{equation}
  for some class $\mathscr K$ function $\phi$ and all times $t\geq 0$.
  Hence for every initial condition and essentially bounded input $u$
  the solution of our system (\ref{eq:ps3}) exists for all times
  $t\geq0$, since $s_\infty$ in (\ref{eq:s_infty}) does not depend on
  $t$.

  \emph{Establishing ISS:} We now utilize an idea from \cite{JTP94}:
  Instead of estimating $|x_i(t)|$ with respect to $|x_i(0)|$ in
  (\ref{eq:ps2}), we can also have the point of view that our
  trajectory started in $x_i(\tau)$ at time $0\leq\tau\leq t$ and we
  followed it for some time $t-\tau$ and reach $x_i(t)$ at time $t$.
  For $\tau=\nicefrac t2$ this reads
  \begin{eqnarray}
    |x_i(t)|&\leq&\beta_i(|x_i(\nicefrac t2)|,\nicefrac t2)+\sum_{j\ne
     i}\gamma_{ij}(\|x_{i,[\nicefrac t2,t]}\|_\infty) +\gamma(u)
   \nonumber
   \\
     &\leq&\beta_i(s_\infty,\nicefrac t2)+\sum_{j\ne
     i}\gamma_{ij}(\|x_{i,[\nicefrac t2,\infty)}\|_\infty) +\gamma(u)
   \label{eq:monotone abschaetzung gemacht}
   \\
   &=&\tilde\beta_i(s_\infty,t)+\sum_{j\ne
     i}\gamma_{ij}(\|x_{i,[\nicefrac t2,\infty)}\|_\infty) +\gamma(u)
   \label{eq:scaled ISS condition}
  \end{eqnarray}
  where we again applied (\ref{eq:monot}) to obtain (\ref{eq:monotone
    abschaetzung gemacht}) and defined
  $$\tilde\beta_i(s_i,t):=\beta_i(s_i,\nicefrac t2),$$
  which is of class $\mathscr {KL}$.

  To write inequality (\ref{eq:scaled ISS condition}) in vector form,
  we define
  \begin{equation}
    \label{eq:definition betatilde}
    \tilde\beta(s,t)
    :=
    \left(
      \begin{array}{c}
        \tilde\beta_1(s_1,t)\\
        \vdots\\
        \tilde\beta_n(s_n,t)
      \end{array}
    \right)
  \end{equation}
  for all $s\in\R^n_+$. Denoting by
  $s^n_\infty:=(s_\infty,\ldots,s_\infty)^T$ we obtain the vector
  formulation of (\ref{eq:scaled ISS condition}) as
  \begin{eqnarray}
    \label{eq:scaled vector ISS}
    \einsnorm{x}{t}\leq \tilde\beta(s^n_\infty,t)+\Gamma\circ
    \supnorm{x}{[\nicefrac t2,\infty)}+\gamma^n(\|u\|_\infty).
  \end{eqnarray}

  \let\meet=\wedge \let\join=\vee
  \let\bigmeet=\bigwedge \let\bigjoin=\bigvee
  \let\epsilon=\varepsilon

  By the boundedness of the solution we can take the upper limit on
  both sides of (\ref{eq:scaled vector ISS}).
  By Lemma~\ref{limsuplemma} we have
  $$\limsup_{t\to\infty} \einsnorm{x}{t}=\limsup_{t\to\infty} \supnorm
  {x}{[\nicefrac t2,\infty)}=: l(x)\,,$$
  and it follows that
  $$(\mbox{Id}-\Gamma)\circ l(x)\leq \gamma^n(\|u\|_\infty)$$
  since
  $\lim_{t\to\infty}\tilde\beta(s^n_\infty,t)=0$. Finally,
  by Lemma~\ref{lemma:boundedness_ID-Gamma} we have
  \begin{equation}
    |l(x)|\leq \phi(|\gamma^n(\|u\|_\infty)|)\label{eq:yields AG}
  \end{equation}
  for some $\phi$ of class $\mathscr K_\infty$.  But (\ref{eq:yields AG}) is the
  asymptotic gain property (\ref{eq:AG}).
  
  Now $0$-GAS is established as follows: First note that for $u\equiv 0$ the
  quantity $s_\infty$ in \eqref{eq:s_infty} is a ${\cal K}$ function of
  $|x(0)|$. So \eqref{eq:s_infty} shows (Lyapunov) stability of the system in
  the case $u\equiv 0$. Furthermore, \eqref{eq:yields AG} shows attractivity
  of $x=0$ for the system \eqref{eq:ps3} in the case $u\equiv 0$. This shows
  global asymptotic stability of $x=0$.

  Hence system (\ref{eq:ps3}) is AG and 0-GAS, which together were proved to be
  equivalent to ISS in \cite[Theorem 1]{SoW96}.
\end{proof}

\subsection{Linear gains and an ISS-Lyapunov version}
\label{subsec:linear-gains-and-iss-lyapunov-function}

Suppose the gain functions $\gamma_{ij}$ are all linear, hence
$\Gamma$ is a linear mapping and (\ref{eq:ps_gamma}) is just
matrix-vector multiplication. Then we have the following
\begin{corollary}
  \label{cor:linear-gains}
  Consider $n$ interconnected ISS systems as in the previous section
  on the problem description with a linear
  
  gain matrix
  $\Gamma$, such that for the spectral radius $\rho$ of $\Gamma$ we have
  \begin{equation}
    \rho(\Gamma)<1.\label{eq:rho<1}
  \end{equation}
 Then the system defined by (\ref{eq:ps3}) is ISS from $u$ to $x$.
\end{corollary}

\begin{remark}
  For non-negative matrices $\Gamma$ 
  
  it is well known that (see, e.g., \cite[Theorem 2.1.1, page
  26, and Theorem 2.1.11, page 28]{BeP79})
  \begin{enumerate}
  \item $\rho(\Gamma)$ is an eigenvalue of $\Gamma$ and $\Gamma$
    possesses a non-negative eigenvector corresponding to
    $\rho(\Gamma)$,
  \item $\alpha x\leq\Gamma x$ holds for some $x\in\R^n_+\setminus\{0\}$
    if and only if $\alpha\leq\rho(\Gamma)$.

  \end{enumerate}
  Hence $\rho(\Gamma)<1$ if and only if $\Gamma s\not\geq s$ for all
  $s\in\R^n_+\setminus\{0\}$.   

  Also, by continuity of the spectrum it is  clear that for such 
  $\Gamma, \rho(\Gamma) <1,$ there always exists a matrix
  $D=\mbox{diag}(1+\alpha_1,\ldots,1+\alpha_n)$ with $\alpha_i>0,
  i=1,\ldots,n$, such that $\Gamma Ds\not\geq s$ for all
  $s\in\R^n_+\setminus\{0\}$. 
\end{remark}

\begin{remark}
  \label{remark:Vidyasagar}
  For the case of large-scale interconnected input-output systems a
  similar result exists, which can be found in a monograph by
  Vidyasagar, cf.~\cite[p. 110]{Vid81}. It also covers Corollary
  \ref{cor:linear-gains} as a special case.  The condition on the
  spectral radius is quite the same, although it is applied to a
  \emph{test matrix}, whose entries are finite gains of products of
  \emph{interconnection operators} and corresponding \emph{subsystem
    operators}. These gains are non-negative numbers and, roughly
  speaking, defined as the minimal possible slope of affine bounds on
  the interconnection operators.
\end{remark}

\begin{proof}(of Corollary~\ref{cor:linear-gains})
  The proof is essentially the same as of
  Theorem~\ref{thm:small-gain_for_networks}, but note that instead of
  Lemma~\ref{lemma:boundedness_ID-Gamma} we now directly have
  existence of
  $$(\mbox{Id}-\Gamma)^{-1}=\mbox{Id}+\Gamma+\Gamma^2+\ldots$$
  since $\rho(\Gamma)<1$ and from the power sum expansion it is
  obvious that $(\mbox{Id}-\Gamma)^{-1}$ is a non-decreasing mapping,
  i.e., for $d_1,d_2\geq0$ we have
  $(\mbox{Id}-\Gamma)^{-1}(d_1+d_2)-(\mbox{Id}-\Gamma)^{-1}(d_1)\geq0$.

  Thus at the two places where Lemma~\ref{lemma:boundedness_ID-Gamma} has been
  used we can simply apply $(\mbox{Id}-\Gamma)^{-1}$ to get the desired
  estimates.
\end{proof}

\paragraph{Construction of an ISS-Lyapunov function}
%auto-ignore
There is another approach to describe the ISS property via so called
ISS-Lyapunov functions, cf. \cite{SoW95}.
\begin{definition}
  A smooth function $V$ is said to be an ISS-Lyapunov function of the
  system \eqref{eq:prelim-system} $\dot x=f(x,u),\,f:\R^{N+L}\rightarrow\R^N$ if
  \begin{enumerate}
  \item $V$ is proper, positive-definite,  that is, there exit
    functions $\psi_1,\psi_2$ of class $\mathscr K_\infty$ such that
    \begin{equation}
      \label{eq:1}
      \psi_1(|x|)\leq V(x)\leq \psi_2(|x|),\quad \forall x\in\R^{n_N};
    \end{equation}
  \item there exists a positive-definite function $\alpha$, a class
  $\mathscr K$-function $\chi$, such that
  \begin{equation}
    \label{eq:3}
    V(x)\geq \chi(|u|) \implies \nabla V(x)f(x,u)\leq-\alpha(|x|).
  \end{equation}
  We call the function $\chi$ the \emph{Lyapunov-gain}.
  \end{enumerate}
\end{definition}

In case of linear Lyapunov-gains a Lyapunov function for the interconnected
system can be constructed, given ISS-Lyapunov functions of the
subsystems. Note, that this time we define the gain matrix $\Gamma$
with respect to the Lyapunov-gains $\gamma_{ij}$.

\begin{figure}[ht]
  \begin{center}
    \epsfig{figure=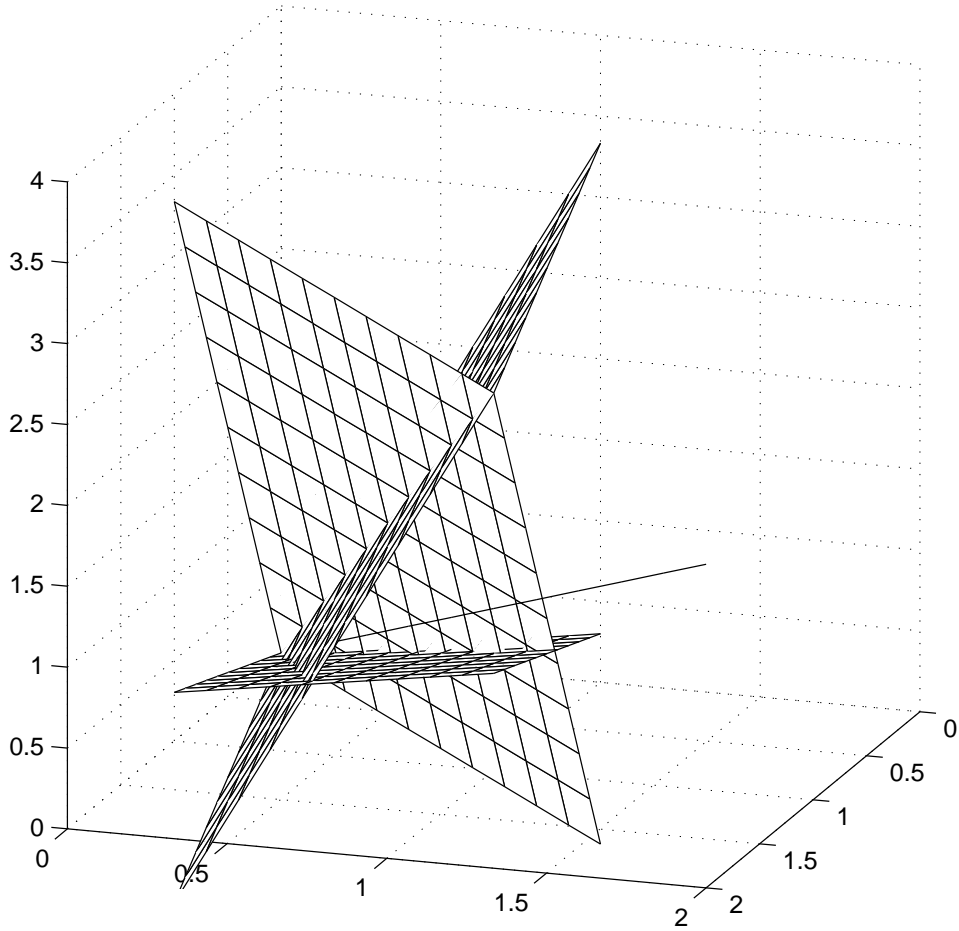, width=.9\textwidth}
  \end{center}
\caption{The sets $\partial\Omega_i$ in $\R^3_+$ and the
eigenvector $s$.}\label{omegas_with_line}
\end{figure}

 Let
$V_1(x_1),\dots,V_n(x_n)$ be some ISS-Lyapunov functions of the
subsystems (\ref{eq:ps1}), allowing for linear Lyapunov-gains
$\gamma_{ij}$, i.e., there are some
$\mathscr K_\infty$ functions $\psi_{i1},\psi_{i2}$
such that
\begin{equation}
\psi_{i1}(|x_i|)\le V_i(x_i)\le \psi_{i2}(|x_i|),\quad
x_i\in\R^{N_i},
\end{equation}
and some positive-definite functions $\alpha_i$ such that
\begin{equation}\label{Lyapunov-functions}
V_i(x_i)>\max \{ \max_j\{\gamma_{ij}V_j(x_j)\}
,\gamma_i(|u|)\}\;\Rightarrow\;\nabla
V_i(x_i)f_i(x_,u)\le-\alpha_i(V_i(x_i)).
\end{equation}
Consider the positive orthant $\R^n_+$, and let $\Omega_i$ be the
subsets of $\R^n_+$ defined by
\begin{equation}
\Omega_i:=\Big\{(v_1,\dots,v_n)\in\R^n_+:\,v_i>\sum_{j=1}^n\gamma_{ij}v_j\Big\}.
\end{equation}
Note that the boundaries $\partial\Omega_i$ are hyperplanes in
case of linear gains.  Now if \eqref{eq:rho<1} or equivalently
$\Gamma s\not\geq s,\,\forall s\in\R^n_+,\, s\ne 0$, holds for
$\Gamma=(\gamma_{ij}),\ i,j=1,\dots,n$, then it follows that
\begin{equation}
\bigcup_{i=1}^n \Omega_i=\mathbb{R}_+^n\setminus \{0\}
    \quad\mbox{and}\quad \bigcap_{i=1}^n \Omega_i\neq\emptyset.
\end{equation}
The proof is the same as of Proposition~\ref{prop:union_omegas},
see below. Thus we may choose an $s> 0$ with $s \in
\bigcap_{i=1}^n \Omega_i,$ which implies that
\begin{equation}\label{property}
s_i>\sum_j \gamma_{ij}s_j,\;i=1,\dots,n;
\end{equation}
see Fig.~\ref{omegas_with_line}. If $\Gamma$ is irreducible, then using
Perron-Frobenius theory we see that we may choose $s$ to be a (positive)
eigenvector of $\Gamma$.

\begin{theorem}
  Let $V_i$ be an ISS-Lyapunov function as in (\ref{Lyapunov-functions}) of the $i^{th}$ subsystem from
  \eqref{eq:ps1}, $i=1,\dots,n,$ and $s$ be a positive vector with \eqref{property}. Then an ISS Lyapunov function of the
  interconnected system \eqref{eq:ps1} is given by
\begin{equation}\label{Lyapunov-function}
V(x_1,\dots,x_n):=\max_{i}{\frac{V_i(x_i)}{s_i}}.
\end{equation}
\end{theorem}
\begin{proof}
Let $\gamma(|u|):=\max_i\gamma_i(|u|)$ which is a $\mathscr K$
class function. In the following we show that there exists a
positive definite function $\alpha$ such that:
\begin{equation}
V(x)\ge\gamma(|u|)\implies\nabla V(x)f(x,u)\le-\alpha(V(x)).
\end{equation}

Let  $M_i$ be open domains in $\R^n_+$ defined by
\begin{equation}\label{M-i}
M_i:=\Big\{(v_1,\dots,v_n)\in\R^n_+:\;\frac{v_i}{s_i}>\max_{j\neq
  i}\Big\{\frac{v_j}{s_j}\Big\}\Big\},
\end{equation}
and  let $P_i$ be the 2-dimensional planes spanned by $s$ and the $i$-th axis, i.e.,
\begin{equation}
P_i=\Big\{v\in\R^n_+\Big\vert\,\frac{v_k}{s_k}=\frac{v_j}{s_j};\;
\forall k,j\neq i\Big\}.
\end{equation}
Note that $V$ defined by (\ref{Lyapunov-function}) is continuous in
$\R^n_+$ and can only fail to be differentiable on the planes $P_i$.

Now take any $\hat x=(\hat x_1,\dots,\hat x_n)\in\R^n$ with
$(V_1(\hat x_1),\dots,V_n(\hat x_n))\in M_i$ then it follows that
in some neighborhood $U$ of $\hat x$  we have
$V(x)=\frac{V_i(x_i)}{s_i}$ for all $x\in U$ and
\begin{equation}
  V_i(x_i)>\max_{j\ne i}\Big\{\frac{s_i}{s_j}V_j(x_j)\Big\} >
  \max_{j\ne i}\{\gamma_{ij}V_j(x_j)\}
\end{equation}
(the last inequality follows from(\ref{property})), hence by
(\ref{Lyapunov-functions}), if $V(x)=V_i(x_i)/s_i>\gamma_i(|u|)$,
then
\begin{equation}
\label{decrease}
  \nabla V(x)f(x,u)=\frac{1}{s_i}\nabla V_i(x_i)f_i(x,u)\le
  -\frac{1}{s_i}\alpha_i(V_i(x_i))<
  -\tilde{\alpha_i}(V(x)),
\end{equation}
where $\tilde{\alpha_i}$ are positive-definite functions, since
$s_i=const>0.$

It remains to consider $x\in\R^n$ such that
$(V_1(x_1),\dots,V_n(x_n))\in \overline{M_i}\cap\overline{M_j}$, where
$V(x)$ may be not differentiable.

For this purpose we use some results from \cite{CLSW98}. For
smooth functions $f_i,\;i=1,\dots,n$ it follows that
$f(x,u)=\max\limits_i \{f_i(x,u)\}$ is Lipschitz and Clarke's
subgradient of $f$ is given by
\begin{equation}
\partial_{Cl} f(x)=co\Big\{\bigcup\limits_{i\in M(x)}\nabla_x f_i(x,u)\Big\},\quad M(x)=\{i:\,f_i(x,u)=f(x)\},
\end{equation}
i.e., in our case
\begin{equation}
\partial_{Cl} V(x)=co\Big\{\frac{1}{s_i}\nabla V_i(x):\;\frac{1}{s_i} V_i(x)=V(x)\Big\}.
\end{equation}
Now for every extremal point of $\partial_{Cl} V(x)$ a decrease
condition is satisfied by \eqref{decrease}. By convexity, the same
is true for every element of $\partial_{Cl} V(x)$. Now
Theorems~4.3.8 and 4.5.5 of \cite{CLSW98} show strong invariance
and attractivity of the set $\{x:\, V(x) \leq \gamma(\|u\|)\}$. It
follows that $V$ is an ISS-Lyapunov function for the
inter\-con\-nection~(\ref{eq:ps1}).
\end{proof}

See also Section~\ref{sec:geom-interpr} for some more considerations
into this directions.

\section{Interpretation of the generalized small-gain condition}
%auto-ignore
\label{sec:interpretations}
In this section we wish to provide insight into the small-gain condition of
Theorem~\ref{thm:small-gain_for_networks}. We first show, that the result
covers the known interconnection results for cascades and feedback
interconnections. We then compare the condition with the linear case.

Further we state some algebraic and graph theoretical relations
and investigate some associated artificial dynamical system
induced by the gain matrix $\Gamma$.  We complete this section
with some geometrical considerations and an overview map of all
these contiguities.

\subsection{Connections to known results}

As an easy consequence of Theorem~\ref{thm:small-gain_for_networks} we
recover, that an arbitrarily long feed forward cascade of ISS subsystems is
ISS again. If the subsystems are enumerated consecutively and the gain
function from subsystem $j$ to subsystem $i>j$ is denoted by $\gamma_{ij}$,
then the resulting gain matrix has non-zero entries only below the diagonal.
For arbitrary $\alpha\in\mathscr K_\infty$ the gain matrix with entries
$\gamma_{ij}\circ (\mbox{Id}_{\R_+}+\alpha)$ for $i>j$ and $0$ for $i\leq j$
clearly satisfies \eqref{eq:Gamma D Id}. Therefore the feed forward cascade
itself is ISS.

Consider $n=2$ in equation (\ref{eq:ps1}), i.e., two subsystems with linear
gains.  Then in Corollary~\ref{cor:linear-gains} we have
$$
\Gamma=\left[\begin{array}{cc}
    0&\gamma_{12}\\ \gamma_{21}&0\end{array}\right],\quad\gamma_{ij}\in\R_+
$$
and $\rho(\Gamma)<1$ if and only if $\gamma_{12}\gamma_{21}<1.$
Hence we obtain the known small-gain theorem, cf. \cite{JMW96} and \cite{Gru02}.

For nonlinear gains and $n=2$ the condition (\ref{eq:Gamma D Id})
in Theorem \ref{thm:small-gain_for_networks} reads as follows: There
exist $\alpha_1, \alpha_2\in\mathscr K_\infty$ such that
$$
\left(\begin{array}{c}
    \gamma_{12}\circ(\mbox{Id}+\alpha_2)(s_2)\\
    \gamma_{21}\circ(\mbox{Id}+\alpha_1)(s_1)
  \end{array}\right)\ngeq
\left(\begin{array}{c} s_1 \\ s_2\end{array}\right),
$$
for all $(s_1,s_2)^T \in\R^2_+$. This is easily seen to be equivalent to
$$\gamma_{12}\circ(\mbox{Id}+\alpha_2)\circ\gamma_{21}\circ(\mbox{Id}+\alpha_1)(s)<s\,,\quad
\forall s>0.$$
To this end it suffices to check what happens for to the vector
$[\gamma_{12}\circ(\mbox{Id}+\alpha_2)(s_2), s_2]^T$ under $\Gamma$ along with
a few similar considerations.  The latter is equivalent to the condition in
the small-gain theorem of \cite{JTP94}, namely, that for some
$\tilde\alpha_1,\tilde\alpha_2\in\mathscr K_\infty$ it should hold that
\begin{equation}
  (\mbox{Id} +
  \tilde\alpha_1) \circ \gamma_{21} \circ (\mbox{Id} + \tilde\alpha_2)
  \circ \gamma_{12}(s) \leq s\,
  ,\quad\forall s>0\,,
  \label{eq:JTP94-Bedingung}
\end{equation}
for all $s\in\R_+$, hence our theorem contains this result as a
particular case.

\begin{example}
  \label{JTP-example}
  The condition (\ref{eq:JTP94-Bedingung}) of \cite{JTP94} seems to be
  very similar to the small-gain condition
  $\gamma_{12}\circ\gamma_{21}(s)<s$ of \cite{JMW96} and \cite{Gru02},
  however those $\gamma$'s have some different meanings in these
  papers. This similarity raises the question, whether the
  compositions with $(\mbox{Id}+\tilde\alpha_i),\;i=1,2$ in
  (\ref{eq:JTP94-Bedingung}) or more generally with $D$ in
  \eqref{eq:Gamma D Id} is necessary. The answer is positive. Namely,
  there is a counterexample providing a system of two ISS subsystems
  with $\gamma_{12}\circ\gamma_{21}(s)<s$ ($\gamma$'s are defined as
  above) which is not ISS.

  Consider the equation
  $$
  \dot x=-x+u(1-e^{-u}),\quad x(0)=x^0\in\R, u\in\R.
  $$
  Integrating it follows
  $$
  x(t)= e^{-t}x^0+\int_0^t
  e^{-(t-\tau)}u(\tau)(1-e^{-u(\tau)})\,d\tau
  $$
  $$
  \le e^{-t}x^0+||u||_\infty(1-e^{-||u||_\infty})=e^{-t}x^0+\gamma(\|u\|_\infty),\quad
  \gamma(s)<s.
  $$
  Then for a feedback system
  \begin{eqnarray}
    \dot x_1&=&-x_1+x_2(1-e^{-x_2})+u(t),\\
    \dot x_2&=&-x_2+x_1(1-e^{-x_1})+u(t)
  \end{eqnarray}
  we have ISS for each subsystem with $x_i(t)\le
  e^{-t}x^0_i+\gamma_i(||x_i||)+\eta_i(||u||),$ where $\gamma_i(s)<s$
  and hence $\gamma_1\circ\gamma_2(s)<s$ for $s>0$, but there is a solution
  $x_1=x_2=const$, i.e.,
  $$
  \dot x_1=-x_2 e^{-x_2}+u, \quad\mbox{with}\; u=x_2 e^{-x_2},
  $$
  and $x_1=x_2$ can be chosen arbitrary large with $u\rightarrow 0$
  for $x_1\rightarrow\infty.$ Hence the condition $\Gamma(s)\not\geq
  s$, for all $s\in\R^n_+\setminus\{0\}$, or for two subsystems
  $\gamma_{12}\circ\gamma_{21}(s)<s$, for all $s>0$, is not
  sufficient for the input-to-state stability of the composite system
  in the nonlinear case.
\end{example}

\subsection{Algebraic Interpretation}
\label{sec:algebr-interpr}

In this subsection we first relate the network small-gain condition
\eqref{eq:Gamma D Id} to well known properties of matrices in the
linear case. This gives some idea how the new condition can be
understood and what subtle differences appear in the nonlinear case.
Then we extend some graph theoretical results for non-negative
matrices to nonlinear gain matrices. These are needed later on.

For a start, we discuss some algebraic consequences from
\eqref{eq:Gamma D Id}. Recall that for a
non-negative
matrix $\Gamma$ the following are equivalent:
\begin{enumerate}
\item 
  $\rho(\Gamma) < 1$,
\item
  $\forall s\in \R^n_+\setminus \{ 0 \} : \Gamma s\not \geq s$,
\item
  $\Gamma^k \to 0$, for $k\to \infty$,
\item
  there exist $a_1,\ldots,a_n>0$ such that $\forall s\in
  \R^n_+\setminus \{ 0 \}$:
  \begin{equation*}
    \Gamma (I +
    \mathrm{diag}(a_1,\ldots,a_n)) s\not \geq s.
  \end{equation*}
\end{enumerate}

Note that (iv) is the linear version of (\ref{eq:Gamma D Id}). As
condition (i) is not useful in the nonlinear setting, we have turned
to (ii), which we later strengthened to \eqref{eq:Gamma D Id}.

In the nonlinear case we find the obvious implication:
\begin{proposition}
  Condition (\ref{eq:Gamma D Id}) implies that
  \begin{equation}
    \label{eq:Interp-1}
    \Gamma(s)\ngeq s \quad\mbox{for any}\quad s\in\R^n_+\setminus\{0\}.
  \end{equation}
\end{proposition}
\begin{proof}
  By the monotonicity of $\Gamma$ it is obvious that (\ref{eq:Gamma D
    Id}) implies \eqref{eq:Interp-1}.
\end{proof}

Note that the contrary is not true:

\begin{example}
  Let
  $$\gamma_{12}=\mbox{Id}_{\R_+}$$
  and $$\gamma_{21}(r)=r(1-e^{-r}).$$
  Since already
  $\lim_{r\to\infty}(\gamma_{12}\circ\gamma_{21}-\mbox{Id})(r)=0$
  there are certainly no class $\mathscr K_\infty$ functions
  $\tilde\alpha_i, i=1,2$ such that (\ref{eq:JTP94-Bedingung}) holds.
\end{example}

\begin{remark} We like to point out the connections between 
  non-negative matrices, our gain matrix and directed graphs.

  A (finite) directed graph $G=\{V,E\}$ consists of a set $V$ of
  vertices and a set of edges $E\subset V\times V$. We may identify
  $V=\{1,\ldots,n\}$ in case of $n$ vertices. The adjacency matrix
  $A_G=(a_{ij})$ of this graph is defined by $$a_{ij}=
  \begin{cases}
    1& \mbox{if}~(i,j)\in E,\\
    0& else.
  \end{cases}$$
  The other way round, given an  $n\times n$-matrix $A$, one defines
  the graph $G(A)=\{V,E\}$ by $V:=\{1,\ldots,n\}$ and $E=\{(i,j)\in V\times
  V: a_{ij}\ne 0\}$.
  
  There are several concepts and results of (non-negative) matrix
  theory, which are of purely graph theoretical nature. Hence the same
  can be done for our interconnection gain matrix $\Gamma$. We may
  associate a graph $G(\Gamma)$, which represents the interconnections
  between the subsystems, in the same manner, as we would do for
  matrices.

  We could also use the graph of the transpose of $\Gamma$ here for
  compatibility with our previous notation ($\gamma_{ij}$ encodes
  whether or not subsystem $j$ influences subsystem $i$) and the
  standard notation in graph theory (edge from $i$ to $j$), then the
  arrows in $G(\Gamma)$ would point in the \lq{}right\rq{}\ direction. But
  this does not affect the following results.

  For instance, we say $\Gamma$ is \emph{irreducible}, if $G(\Gamma)$ is
  \emph{strongly connected}, that is, for every pair of vertices
  $(i,j)$ there exists a sequence of edges (a \emph{path})
  connecting vertex $i$ to vertex $j$. Obviously $\Gamma$ is
  irreducible if and only if $\Gamma^T$ is. $\Gamma$ is called
  \emph{reducible} if it is not irreducible.

  The gain matrix $\Gamma$ is \emph{primitive}, if its associated
  graph $G_\Gamma:=G(\Gamma)$ is primitive, i.e., there exists a
  positive integer $m$ such that $(A_{G_\Gamma})^m$ has only positive
  entries.

  These definitions and the following important facts can be found in
  \cite{BeP79} and only depend on the associated graph.

  If $\Gamma$ is reducible, then a permutation transforms it into a block
  upper triangular matrix. From an interconnection point of view, this splits
  the system into cascades of subsystems each with irreducible adjacency
  matrix.

\end{remark}

\begin{lemma}
  \renewcommand{\labelenumi}{\alph{enumi})}
  \label{lemma:BeP79-adjacenzmatrix resultate}
  Assume the gain matrix $\Gamma$ is irreducible. Then there are two
  distinct cases:
  \begin{enumerate}
  \item 
    The gain matrix
    $\Gamma=(\gamma_{ij}(\cdot)),$ where $\gamma_{ij}(\cdot)\in\mathscr{K}$
    or $\gamma_{ij}=0$, is primitive and hence there is a
    non-negative integer $k_0$ such
    that $\Gamma^{k_0}$ has elements
    $\gamma^{k_0}_{ij}(\cdot)\in\mathscr{K}$ for any $i,j$.
  \item 
    The gain matrix $\Gamma$ can be transformed to
    \begin{equation}
      P\Gamma P^T=\left(\begin{array}{ccccc}
          0     & A_{12} &  0     & \dots   & 0\\
          0     &0       & A_{23} & \dots   & 0\\
          \vdots&        &\vdots &\ddots&\vdots\\
          0     &0       &0      &\dots&A_{\nu-1,\nu}\\
          A_{\nu1}&0       &0      &\dots&0
        \end{array}\right)\label{eq:cogredience form}
    \end{equation}
    using some permutation matrix $P$, where the zero blocks on the
    diagonal are square and where $\Gamma^{\nu}$ is of block
    diagonal form with square primitive blocks  on the
    diagonal.
  \end{enumerate}
\end{lemma}

\begin{proof}
  Let $A_{G_\Gamma}$ be the adjacency matrix corresponding to the
  graph associated with $\Gamma$. This matrix is primitive if and only
  if $\Gamma$ is primitive.  Note that the $(i,j)^{\mbox{\small{th}}}$
  entry of $A_{G_\Gamma}^k$ is zero if and only if the
  $(i,j)^{\mbox{\small{th}}}$ entry of $\Gamma^k$ is zero.
  Multiplication of $\Gamma$ by a permutation matrix only rearranges
  the positions of the class $\mathscr K$-functions, hence this
  operation is well defined. From these considerations it is clear,
  that it is sufficient to prove the lemma for the matrix
  $A:=A_{G_\Gamma}$. But for non-negative matrices this result is an
  aggregation of known facts from the theory of non-negative matrices,
  see, e.g., \cite{BeP79} or \cite{LaT85}.
\end{proof}
\subsection{Asymptotic Behavior of $\Gamma^k$}

A related question to the stability of the composite system
\eqref{eq:ps3} is, whether or not the discrete positive dynamical
system defined by
\begin{equation}
  \label{eq:int:dynsys s=Gs}
  s_{k+1}=\Gamma(s_k),\quad k=1,2,\ldots
\end{equation}
with given initial state $s_0\in\R^n_+$ is globally asymptotically
stable. Under the assumptions we made for
Theorem~\ref{thm:small-gain_for_networks} this is indeed true for irreducible $\Gamma$.

\begin{theorem}
  \label{thm:ass-dynamical-system}
  Assume that $\Gamma$ is irreducible. Then the system defined by
  \eqref{eq:int:dynsys s=Gs} is globally asymptotically stable if and
  only if $\Gamma(s)\not\geq s$ for all $s\in\R^n_+\setminus\{0\}$.
\end{theorem}
The proof will make use of the following result:

\begin{proposition}
  \label{prop:uniformly_attractive}
  The condition
  \begin{equation}
    \label{eq:Interp-2}
    \lim\limits_{k\rightarrow\infty}\Gamma^k(s)\rightarrow0 \quad
    \mbox{for any fixed}\quad s\in\R^n_+
  \end{equation}
  implies \eqref{eq:Interp-1}. Moreover if $\Gamma$ is irreducible,
  then both are equivalent.
\end{proposition}

Note that the converse implication is generally not true for reducible
maps $\Gamma$, such that \eqref{eq:Interp-1} holds. See
Example~\ref{example:counterxample_gamma_not_to_null}. But it is
trivially true, if $\Gamma$ is linear.

\begin{proof}
  Condition \eqref{eq:Interp-1} follows from \eqref{eq:Interp-2},
  since if $\Gamma(s_0)\ge s_0$ for some $s_0\in\R^n_+\setminus\{0\}$
  then $\Gamma^k(s_0)\ge \Gamma^{k-1}(s_0)\geq s_0$ for
  $k=2,3,\ldots$. Hence the sequence $\{\Gamma^k(s_0)\}_{k=0}^\infty$
  does not converge to $0$.

  Conversely, assume that \eqref{eq:Interp-1} holds and that $\Gamma$
  is irreducible.\\
  Step 1. First we prove that
  for any $s\in\R^n_+\setminus\{0\}$
  \begin{equation}\label{eg:gamma}
    \Gamma^k(s)\ngeq s,\quad k\in\mathbb{N}.
  \end{equation}
  Assume there exist some $k>1$ and $s\neq0$ with $\Gamma^k(s)\geq s.$ Define
  $z\in\R_+^n$ as
  $$
  z:=\max_{l=0,\dots,k-1}\{\Gamma^l (s)\} \,
  \genfrac{}{}{0pt}{}{\geq}{\ne} \, 0\,.
  $$
  By \eqref{eq:monot} and using $\Gamma^k(s)\geq s$  we have
  $$\Gamma(z)\ge\max_{l=1,\dots,k}\{\Gamma^l s\}=\max_{l=0,\dots,k}\{\Gamma^l s\}\ge
  \max_{l=0,\dots,k-1}\{\Gamma^l s\}=z.$$
  This contradicts~\eqref{eq:Interp-1}.

  Step 2. For any fixed $s$ we prove that
  $\limsup_{k\rightarrow\infty}|\Gamma^k(s)| <\infty$.  By
  Lemma~\ref{lemma:BeP79-adjacenzmatrix resultate} we have two cases.
  We only consider case a), then case b) follows with a slight
  modification.

  Assume that $s$ is such that $\limsup_{k\rightarrow\infty}
  |\Gamma^k(s)| = \infty$.  For $t_i >0$ denote the $i^{\mbox{th}}$
  column of $\Gamma^{k_0}$ by
  $$
  \Gamma_i^{k_0}(t_i)=\left(\begin{array}{c}
      \gamma_{1i}^{k_0}(t_i)\\
      \vdots \\
      \gamma_{ni}^{k_0}(t_i)
    \end{array}\right)
  $$

  As $\Gamma^{k_0}$ has no zero entries, for $i=1,\dots,n$ there are
  $T_i\in\R_+$ such that
  \begin{equation}\label{eq:Gko>s}
    \Gamma_i^{k_0}(t_i)>s\quad\mbox{for any}\quad t_i>T_i.
  \end{equation}
  If $|\Gamma^k(s)| \to \infty$ there exists a $k_1$ and an index
  $i_1$ such that
  \[ \Gamma^{k_1}(s)_{i_1} \geq T_{i_1} \,.\]
  The vector $\Gamma^{k_0}\circ\Gamma^k(s),$ seen as a sum of columns, is
  greater than the maximum over these columns, i.e.,
  \begin{equation}\label{eq:max-of-columns}
    \begin{split}
      \Gamma^{k_0}\circ\Gamma^{k_1}(s)\ge\max_i
      \Gamma^{k_0}_i(\Gamma^{k_1}(s)_i)
      \geq \Gamma^{k_0}_{i_1}(\Gamma^{k_1}(s)_{i_1})\\ \geq
      \Gamma^{k_0}_{i_1}(T_{i_1}) \geq s .
    \end{split}
  \end{equation}
  This contradicts Step~1.

  Step 3. So $\left\{\Gamma^k(s)\right\}_{k\ge 1}$ is bounded for any
  fixed $s\in\R^n_+$.  The omega-limit set $\omega(s)$ is defined by
  \begin{equation}\label{omega-lim}
    \begin{split}
      \omega(s)=\Big\{x\;\big\vert\;\exists~\mbox{subsequence}~
      \{k_j\}_{j=1,2,\ldots} \\
      \mbox{such that}~
      \Gamma^{k_j}(s)\xrightarrow{j\rightarrow\infty} x\Big\}.
    \end{split}
  \end{equation}

  This set is not empty by boundedness of
  $\left\{\Gamma^k(s)\right\}_{k\ge 1}$.  The following properties
  follow from this definition and boundedness of the set
  $\left\{\Gamma^k(s)\right\}_{k\ge 1}:$
  $$\forall\, x\in\omega(s)\;\Rightarrow\; \Gamma(x)\in\omega(s),$$
  $$\forall\, x\in\omega(s)\;\exists\, y\in\omega(s):\;\Gamma(y)=x.$$
  I.e., $\omega(s)$ is invariant under $\Gamma.$ The boundedness of $\omega(s)$
  allows to define a finite vector
  $$z=\sup\omega(s).$$
  Then for any $x\in\omega(s)$ it follows $z\ge
  x$ and hence $\Gamma(z)\ge\Gamma(x)$.  Let $y\in\omega(s)$ be such
  that $\Gamma(x)=y$. Then $\Gamma(z)\ge y$.  By the invariance of
  $\omega(s)$ it follows that
  $$\Gamma(z)\ge\sup\{\Gamma(x)\,\vert\,\forall\,x\in\omega(s)\}=z.$$
  This contradicts \eqref{eq:Interp-1} if $z\neq 0$, i.e., $\omega(s)=\{0\}$. This is
  true for any $s\in\R^n_+$. Hence \eqref{eq:Interp-2} is proved as a consequence of
  \eqref{eq:Interp-1}, provided that $\Gamma$ is irreducible.
\end{proof}

\begin{example}
  \label{example:counterxample_gamma_not_to_null}
Consider the map $\Gamma:\R^2_+\to\R^2_+$ defined by
$$\Gamma:=
\begin{bmatrix}
  \gamma_{11}&\mbox{id}\\
  0&\gamma_{22}
\end{bmatrix}
$$ where for $t\in\R_+$ $$\gamma_{11}(t):=t(1-e^{-t})$$
and the function $\gamma_{22}$ is constructed in the sequel.  First
note that $\gamma_{11}\in\mathcal{K}_\infty$ and
$\gamma_{11}(t)<t,\,\forall t>0$.  Let $\{\epsilon_k\}_{k=1}^\infty$ a strictly
decreasing sequence of positive real numbers, such that
$\lim_{k\to\infty}\epsilon_k=0$ and
$\lim_{K\to\infty}\sum_{k=1}^K\epsilon_k=\infty$.  For $k=1,2,\ldots$
define
$$\gamma_{22}\bigg( \epsilon_k + (1+\sum_{j=1}^{k-1}\epsilon_j)
e^{-(1+\sum_{j=1}^{k-1}\epsilon_j)}\bigg) :=
\epsilon_{k+1} + (1+\sum_{j=1}^{k}\epsilon_j)
e^{-(1+\sum_{j=1}^{k}\epsilon_j)}$$ and observe that
$$ \epsilon_k + (1+\sum_{j=1}^{k-1}\epsilon_j)
e^{-(1+\sum_{j=1}^{k-1}\epsilon_j)} >
\epsilon_{k+1} + (1+\sum_{j=1}^{k}\epsilon_j)
e^{-(1+\sum_{j=1}^{k}\epsilon_j)}\,,$$
since $\epsilon_k>\epsilon_{k+1}$ for all $k=1,2,\ldots$ and the map $t\mapsto
t\cdot e^{-t}$ is strictly decreasing on $(1,\infty)$.

Moreover we have by assumption, that $$ \epsilon_k +
(1+\sum_{j=1}^{k-1}\epsilon_j)
e^{-(1+\sum_{j=1}^{k-1}\epsilon_j)}\xrightarrow[k\to\infty]{}0.$$

These facts together imply that $\gamma_{22}$ may be extrapolated to some
$\mathcal{K}_\infty$-function, in a way such that
$\gamma_{22}(t)<t,\,\forall t>0$ holds.

Note that by our particular construction we have $\Gamma(s)\not\geq s$
for all $s\in\R^2_+\setminus\{0\}$.
Now define $s^1\in\R^2_+$  by $$s^1:=
\begin{bmatrix}
  1\\
  1+e^{-1}
\end{bmatrix}$$ and for $k=1,2,\ldots$ inductively define
$s^{k+1}:=\Gamma(s^k)\in\R^2_+$.

By induction one verifies that $$s^{k+1}=\Gamma^k(s^1)=
\begin{bmatrix}
  1+\sum_{j=1}^{k}\epsilon_j\\
  \epsilon_{k+1} + (1+\sum_{j=1}^{k}\epsilon_j)
  e^{-(1+\sum_{j=1}^{k}\epsilon_j)}
\end{bmatrix}.$$

By our previous considerations and assumptions we easily obtain that
the second component of the sequence $\{s^k\}_{k=1}^\infty$ strictly
decreases and converges to zero as $k$ tends to infinity. But at the
same time the first component strictly increases above any given
bound.

Hence we established that $\Gamma(s)\not\geq s~\forall s\ne0$ in
general does not imply $\forall s\ne 0:~\Gamma^k(s)\to0$ as $k\to\infty$.
\end{example}

\begin{remark}
  Note that we can even turn the constructed 2x2-$\Gamma$ into the
  null-diagonal form, that is assumed in
  Theorem~\ref{thm:small-gain_for_networks}.
  Using the same notation for $\gamma_{ij}$ as in
  Example~\ref{example:counterxample_gamma_not_to_null}, we just
  define
  \[
  \Gamma:=
  \begin{bmatrix}
    0&\gamma_{11}&\mbox{id}&0\\
    \gamma_{11}&0&0&\mbox{id}\\
    0&0&0&\gamma_{22}\\
    0&0&\gamma_{22}&0
  \end{bmatrix}
  \qquad\mbox{and}
  \qquad
  s^1:=
  \begin{bmatrix}
    1\\
    1\\
    1+e^{-1}\\
    1+e^{-1}
  \end{bmatrix}
  \]
  and easily verify that $\Gamma^k(s^1)$ does not converge to $0$.
\end{remark}

\begin{proof}[Proof of Theorem~\ref{thm:ass-dynamical-system}]
  If \eqref{eq:int:dynsys s=Gs} is asymptotically stable, it is in
  particular attracted to zero, so by
  Proposition~\ref{prop:uniformly_attractive} and irreducibility of
  $\Gamma$ we establish \eqref{eq:Interp-1}.

  Conversely assume \eqref{eq:Interp-1}. Clearly $0\in\R^n_+$ is an
  equilibrium point for \eqref{eq:int:dynsys s=Gs} and by
  Proposition~\ref{prop:uniformly_attractive} it is globally
  attractive. It remains to prove stability, i.e., for any
  $\epsilon>0$ there exists a $\delta>0$ such that $|s_0| < \delta$
  implies $\Gamma^k(s_0)<\epsilon$ for all times $k=0,1,2,\ldots$

  Given $\epsilon>0$ we can choose an $r\in\bigcap_{i=1}^n\Omega_i\cap
  S_\epsilon$ where $S_\epsilon$ is the sphere around $0$ of radius
  $\epsilon$ in $\R^n_+$. Define $\delta$ by
  \[ \delta := \sup \{ d \in \R_+: s < r~ \forall s\in B_d(0)\}. \]
  Here $B_d(0)$ denotes the ball of radius less than $d$ in
  $\R^n_+$ around the origin with respect to the Euclidean norm.
  Clearly we have $r>s_0$ for all $|s_0|<\delta$.  Since
  $r\in\bigcap_{i=1}^n\Omega_i\ne\emptyset$ we have $r>\Gamma(r)$ and
  therefore $r>\Gamma(r)\geq\Gamma^2(r)\geq \ldots$ and even
  $\Gamma^k(r)\xrightarrow{k\to\infty}0$ again by
  Proposition~\ref{prop:uniformly_attractive}.

  Hence for any $s_0$ such that $|s_0|<\delta$ we have
  $\Gamma^k(r)\geq\Gamma^k(s_0)$ for all $k=0,1,2,\ldots$ by
  monotonicity of $\Gamma$ and therefrom $\Gamma^k(s_0)<\epsilon$ for
  all $k=0,1,2,\ldots$
\end{proof}

\subsection{Geometrical Interpretation}
\label{sec:geom-interpr}

For the following statement let us define the open domains
$$\Omega_i=\left\{x\in\mathbb{R}^N:\; |x_i| >\sum\limits_{j\neq
    i}\gamma_{ij}(|x_j|)\right\},$$
where $N=\sum_{j=1}^nN_j$ and $x$ is partitioned to
$(x_1,\ldots,x_n)$ with $x_i\in\R^{N_i}, i=1,\ldots,n$, as in
\eqref{eq:ps1}.

\begin{proposition}
  \label{prop:union_omegas}
  Condition \eqref{eq:Interp-1} is equivalent to
  \begin{equation}
    \label{eq:omegas}
    \bigcup_{i=1}^n \Omega_i=\mathbb{R}^N\setminus \{0\}
    \quad \mbox{and}\quad \bigcap_{i=1}^n \Omega_i\neq\emptyset.
  \end{equation}
\end{proposition}

\begin{proof}
    Let $s\neq 0$. Formula \eqref{eq:Interp-1} is equivalent to the
    existence of at least one index $i\in\{1,\dots,n\}$ with
    $s_i>\sum_{j\neq i}\gamma_{ij}(s_j).$ This proves the first part
    of \eqref{eq:omegas}.

    It remains to show, that \eqref{eq:Interp-1} implies
    $\bigcap_{i=1}^n \Omega_i\neq\emptyset$. We may restrict ourselves
    to the positive orthant in $\R^n$, and the sets
    $$\tilde\Omega_i=\left\{ s\in\mathbb{R}^n_+:\; s_i >\sum\limits_{j\neq
        i}\gamma_{ij}(s_j)\right\}$$ instead of $\Omega_i,\ i=1,\ldots,n$.

    For an index set $I$ we define
    $E_I=\{s\in\R^n_+:s_m=0\quad\mbox{for}\quad m\notin I\}$. Note
    that points of $E_I$ can not be in $\tilde\Omega_m$ for $m\notin
    I$.  Consider $\tilde\Omega_i$ and $\tilde\Omega_j$ for any $i\neq
    j$.  The intersections $\tilde\Omega_i\cap E_{\{i,j\}}$ and
    $\tilde\Omega_j\cap E_{\{i,j\}}$ of this two domains with the
    plane $E_{\{i,j\}}$ are nonempty. The points of $\partial\tilde\Omega_i$ lying
    in this plane do not belong to $\tilde\Omega_k$ for any $k\neq j$,
    hence they are in $\tilde\Omega_j$. Since the domains are open it
    follows that the intersections
    $\tilde\Omega_i\cap\tilde\Omega_j\neq\emptyset$ for any $i\neq j.$
    Denote $\tilde\Omega_{ij}=\tilde\Omega_i\cap\tilde\Omega_j$ which
    has nonempty intersection with $E_{\{i,j\}}$ by construction. Take
    any $k\neq i,j$.  Consider $E_{\{i,j,k\}}\supset E_{\{i,j\}}$ which has
    non-empty intersection with $\tilde\Omega_{ij}$. Let $x\in\tilde\Omega_{ij}\cap E_{\{i,j,k\}}.$
    There is some $y\in E_{\{i,j,k\}}$ with $y\notin\Omega_{ij}$
    (say $y\in E_{\{k\}}$). Since $E_{\{i,j,k\}}$ is convex the
    segment $\overline{xy}\subset E_{\{i,j,k\}},$ hence there is some point
    $z\in E_{\{i,j,k\}}$ belonging to $\partial\tilde\Omega_{ij}$, i.e.,
    $E_{\{i,j,k\}}\cap\partial\tilde\Omega_{ij}$ is non-empty.

    The points of $\partial\tilde\Omega_{ij}$, which are not in
    $\tilde\Omega_i,\,\tilde\Omega_j$ and lying in $E_{\{i,j,k\}}$ can
    not belong to $\tilde\Omega_{\nu}, \nu\neq k$.  Hence they are in
    $\tilde\Omega_k$ and it follows
    $\tilde\Omega_i\cap\tilde\Omega_j\cap\tilde\Omega_k\neq\emptyset.$
    By iteration the second part of \eqref{eq:omegas} follows.  
\end{proof}

\begin{figure}[h]
  \centering
  \setlength{\unitlength}{0.07\columnwidth}
  \begin{picture}(10,8)
    \put(0,1){$\displaystyle\Omega_3$}
    \put(4,7.5){$\displaystyle\Omega_2$}
    \put(9.5,3){$\displaystyle\Omega_1$}
    \put(0,0){\epsfig{figure=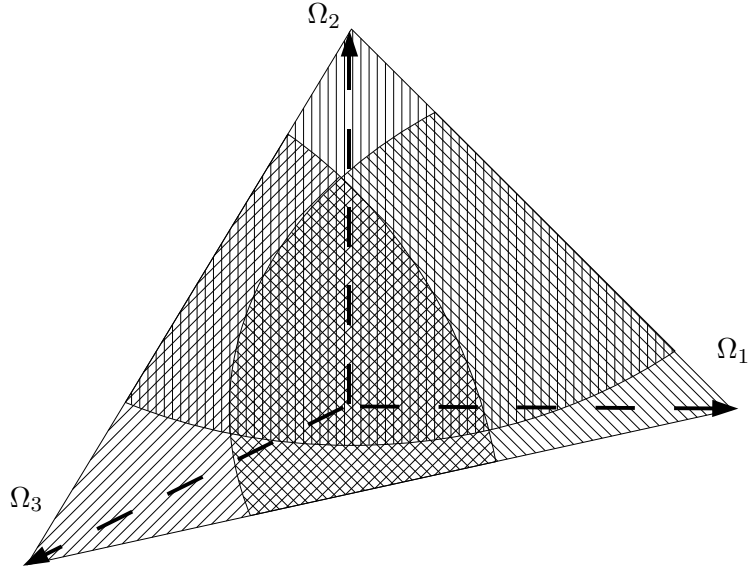, width=.7\columnwidth}}
  \end{picture}
  \caption{Overlapping of
    $\Omega_i$ domains in $\mathbb{R}^3$}
  \label{omegas}
\end{figure}
\setlength{\unitlength}{5mm} 

Let us briefly explain, why the overlapping condition
\eqref{eq:omegas} is interesting: From the theory of ISS-Lyapunov functions it
is known, that a system of the form \eqref{eq:prelim-system} is ISS if and only if
there exists a smooth Lyapunov function $V$ with the property
\[ |x| \geq \gamma(|u|) \Rightarrow \nabla
 V(x) f(x,u) < -W(|x|) \,,\]
for some $W \in {\cal K}$. In the case of our interconnected system this
condition translates to the existence of Lyapunov functions $V_i$ for the
subsystems $i=1,\ldots,n$ with the property
\begin{equation}
\label{Lycond}
  \begin{split}
    |x_i | \geq \sum \gamma_{ij}(|x_j|) + \gamma(|u|)\\
\Rightarrow \nabla V_i(x_i) f_i(x,u) < -W_i(|x_i|) \,,
  \end{split}
\end{equation}
Now for $u=0$ the condition of \eqref{Lycond} is simply, that $x\in \Omega_i$.
Thus the overlapping condition states that in each point of the state space
one of the Lyapunov functions of the subsystems is decreasing. It is an
interesting problem if via this an ISS-Lyapunov function for the whole system may
be constructed.

A typical situation in case of three one dimensional systems ($\mathbb{R}^3$)
is presented on the Figure \ref{omegas} on a plane crossing the positive semi
axis.  The three sectors are the intersections of the $\Omega_i$ with this
plane.

\subsection{Summary map of the interpretations concerning $\Gamma$}

\begin{figure}[th]
\centering
  \[\fbox{
    $ \begin{array}{ccccc}
      &&\exists D~\mbox{as in \eqref{eq:definition von D}}:
      \Gamma\circ D(s)\not\geq s&&\\[1ex]
      &&\Downarrow\ (\Uparrow
      \begin{picture}(0,0)
        \put(0,0){\mbox{~if $\Gamma$ is linear)}}
      \end{picture}
      &&\\[2ex]
      \Gamma^k(s)\xrightarrow{k\to\infty} 0
      &\genfrac{}{}{0pt}{0}{\Rightarrow}{\Leftarrow^*}&\Gamma(s)\not\geq
      s&
      \iff
      &
      \bigcup_{i=1}^n \Omega_i=\R^N\setminus\{0\}
      \\[2ex]
      &&\Updownarrow
      \begin{picture}(0,0)
        \put(0,0){
          \mbox{~if $\Gamma$ is linear}
        }
      \end{picture}&&\\[1ex]
      &&\rho(\Gamma)<1&&
    \end{array}
$}
  \]
  \caption{Some implications and equivalences of the generalized small-gain
  condition.  All statements are supposed to hold for all
    $s\in\R^n_+$, $s\ne0$. The implication denoted by ${}^*$ holds if
    $\Gamma$ is linear or irreducible.}
  \label{fig:implication-diagram}
\end{figure}
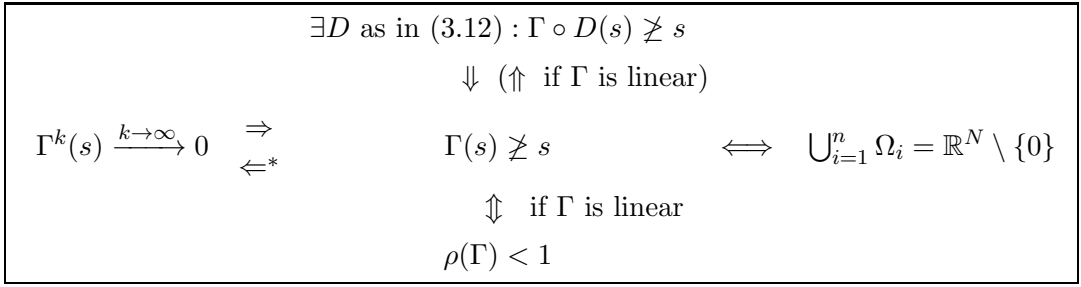

In Figure~\ref{fig:implication-diagram} we summarize the relations
between various statements about $\Gamma$ that were proved in
section~\ref{sec:interpretations}.

\section{Application to linear systems}
\label{sec:appl-line-syst}
%auto-ignore
An important special case is, of course, when the underlying systems are
linear themselves. Consider the following setup, where in the
sequel we omit the external input, formerly denoted by $u$, for
notational simplicity.
Let
\begin{equation}
  \label{eq:als_stable_systems}
  \dot x_j=A_jx_j, \quad x_j\in\R^{N_j}, \quad j=1,\ldots,n
\end{equation}
describe $n$ globally asymptotically stable linear systems, which are
interconnected by the formula
\begin{equation}
  \label{eq:als_interconnected_systems}
  \dot x_j=A_jx_j+\sum_{k=1}^n\Delta_{jk}x_k\quad j=1,\ldots,n,
\end{equation}
which can be rewritten as
\begin{equation}
  \label{eq:als1}
  \dot x=(A +\Delta)x,
\end{equation}
where $A$ is block diagonal, $A=\mbox{diag}(A_j,j=1,\ldots,n)$, each
$A_j$ is Hurwitz (i.e., the spectrum of $A_j$ is contained in the
open left  half plane) and the
matrix $\Delta=(\Delta_{jk})$ is also in block form and encodes the
connections between the $n$ subsystems. We suppose that $\Delta_{jj}=0$
for all $j$. Define the matrix $R=(r_{jk}), R\in\R_+^{n\times n}$, by
$r_{jk}:=||\Delta_{jk}||$.  For each subsystem, there exist
positive constants $M_j,\lambda_j$, such that $e^{A_jt}\leq M_je^{-\lambda_j
t}$ for all $t\geq0$. 

Define a matrix $D\in\R_+^{n\times n}$ by
$D:=\mbox{diag}(\frac{M_j}{\lambda_j},j=1,\ldots,n)$.\\
From the last subsection we obtain
\begin{corollary}
  \label{corollary:anwendung-lineare-systeme}
  If $\rho(D\cdot R)<1$ then (\ref{eq:als1}) is globally
  asymptotically stable.
\end{corollary}

Note that this is a special case of a theorem, which can be
found in Vidyasagar \cite[p. 110]{Vid81}, see Remark
\ref{remark:Vidyasagar}.

\begin{proof}
  Denote the initial value by $x^0$. Then by elementary ODE theory we
  have
  \begin{equation}
    \label{eq:ODE-solution linear subsystem}
    x_j(t)=e^{A_jt}x^0_j+\sum_{k\neq j}
    \int_0^t
    e^{A_j(t-s)} \Delta_{jk}  x_{k}(s) 
    ds
  \end{equation}
  and by standard estimates
  \begin{equation}\label{eq:estimate}
    \vert x_j(t) \vert \le M_j e^{-\lambda_j t} + \sum_{k\neq j}r_{jk}
    \frac{M_k}{\lambda_k} ||x_{k,[0,t]}||.
  \end{equation}
  As one can see from (\ref{eq:estimate}), in this case the gain
  matrix happens to be $\Gamma=D\cdot R$.  
\end{proof}

It is noteworthy, that this particular corollary is also a consequence of more
general and precise results of a recent paper \cite{KHP05} by Hinrichsen,
Karow and Pritchard.

\section{Conclusions}
\label{sec:conclusions}
%auto-ignore
We considered a composite system consisting of an arbitrary number of
nonlinear arbitrarily interconnected subsystems, as they arise in
applications.

For this general case we derived a multisystem version of the
\emph{nonlinear small-gain theorem}.  For
the special case of linear interconnection gains this is a special
case of a known Theorem, cf. \cite[page 110]{Vid81}.
We also showed how our generalized small-gain theorem for networks can
be applied to linear systems.

Many interesting questions remain, for instance concerning the construction of
Lyapunov functions in case of nonlinear (Lyapunov-)gain functions.

%auto-ignore
\section{Acknowledgements}

This research is funded by the German Research Foundation (DFG) as
part of the Collaborative Research Centre 637 "Autonomous Cooperating
Logistic Processes: A Paradigm Shift and its Limitations" (SFB 637).

\bibliographystyle{plain}
\bibliography{literatur}

\end{document}